\documentclass[preprints,article,accept,moreauthors,pdftex]{my-mdpi}

\usepackage[labelformat=simple]{subfig}

\firstpage{1} 
\pubvolume{10}
\issuenum{3}
\articlenumber{135}
\pubyear{2021}
\copyrightyear{2021}
\externaleditor{Academic Editor: {Ioannis Dassios}} 
\datereceived{7 June 2021} 
\daterevised{22 June 2021} 
\dateaccepted{25 June 2021} 
\datepublished{27 June 2021} 
\hreflink{https://doi.org/10.3390/\newline axioms10030135} 
\doinum{10.3390/axioms10030135}

\pdfoutput=1

\Title{Mathematical Analysis of a Fractional COVID-19 Model Applied to Wuhan, Spain and Portugal}

\TitleCitation{Mathematical Analysis of a Fractional COVID-19 Model Applied to Wuhan, Spain and Portugal}

\Author{Fa\"{\i}\c{c}al Nda\"{\i}rou $^{\dagger}$\orcidA{}
and  Delfim F. M. Torres *\orcidB{}}

\AuthorNames{Fa\"{\i}\c{c}al Nda\"{\i}rou and  Delfim F. M. Torres}

\AuthorCitation{Nda\"{\i}rou, F.;  Torres, D.~F.~M.}

\address[1]{Center for Research and Development in Mathematics and Applications (CIDMA), 
Department of Mathematics, University of Aveiro, 3810-193 Aveiro, Portugal}

\corres{\hangafter=1 \hangindent=1.05em \hspace{-0.82em}Correspondence: delfim@ua.pt; Tel.: +351-234-370-668}

\firstnote{\hangafter=1 \hangindent=1.05em \hspace{-0.82em}This research is part of first author's Ph.D. project, 
which is carried out at the University of Aveiro under 
the Doctoral Program in Applied Mathematics
of Universities of Minho, Aveiro, and Porto (MAP-PDMA).} 

\abstract{We propose a qualitative analysis of a recent fractional-order COVID-19 model. 
We start by showing that the model is mathematically and biologically well posed. 
Then, we give a proof on the global stability of the disease free equilibrium point. 
Finally, some numerical simulations are performed to ensure stability  
and convergence of the disease free equilibrium point.}



\begin{document}


\section{Introduction}

The study and analysis of the spread of infectious diseases via mathematical models has become an
important tool to understand diseases' epidemiological prototypes, such as prevalence, duration 
of the epidemic, and its burden on a population. A variety of techniques can be considered in the 
modeling formulation in order to reflect particularities in the transmission routes of a disease. 
In general, there exist two major routes:
\begin{enumerate}
\item[(i)] direct contact, that is, from infected to non-infected individuals 
(for instance, through body fluids for Ebola virus disease \cite{atangana2014, Ngwa, MR3349757} 
or through sexual intercourse for HIV infection \cite{MR3808514,OC:HIV:PLoSCompBio:2015,SilvaTorres:TBHIV:2015} 
or by non-linearities in the transmission incidence \cite{HU201212,PMID:2061695});

\item[(ii)] indirect contact, which is due to the presence of carriers in the environment, 
for example, through mosquitoes and their aquatic phase \cite{Moreno:2017,fays4,fays3,MR3557143}. 
\end{enumerate} 

These features are fundamental in the transmission process of any disease in leading to specific 
health status of each individual in the community. Thus, individuals in the community can be classified 
into disjoints compartments according to their specific health~status.

The so-called classical SIR models refer to three epidemiological status (susceptible--infected--recovered) 
and were firstly introduced by Kermack and McKendrick \cite{Kendrick1} and constitute the basis foundation 
of compartmental models. Additionally, some earlier development aiming to introduce demographic factors in compartmental 
models through birth and death rate were considered by these authors \cite{Kendrick2, Kendrick3}. Nowadays, 
mathematical models are more complex and more realistic by including a large number of health status (see, e.g.,~\cite{Area:in:press,MR3808514,Ngwa} for models with much more epidemiological states). However, in relation 
to the analysis of these models, one of the challenges is how to establish stability 
of equilibrium points~\cite{smith, melnikkk, steven, teschl} 
in order to acquire insights into the dynamics behavior of such models. 

Recently, great considerations have been made to models described by fractional differential equations 
in the field of mathematical epidemiology \cite{sulami,MR3999702,fays2}. The most essential property 
of these models is their memory effect \cite{memory}, which is not appearing in the traditional 
instantaneous differential equations. This effect is fully captured by flexibility of the order 
of differentiation for fractional derivatives and can be seen as a hereditary property or a variety 
on strains and genomes of viruses (as conjectured in \cite{fays2}), which is useful for epidemic models. 
In \cite{math9121321}, a model for controlling the Coronavirus pandemic 2019 (COVID-19) in India with 
Caputo--Fabrizio fractional derivatives and the homotopy analysis transform method is investigated.
Often, fractional-order models give rise to theoretical models that allow a significant improvement 
in the fitting of real data, when compared with analogous classical models \cite{MR3872489}.
However, there are also cases where the fractional-order models do not bring any advantage 
\cite{Boukhouima2021}. In the case of the model under investigation, 
the fractional-order model has clear advantages to describe the spread of COVID-19 in Galicia, 
Spain, and Portugal, but does not offer advantages with respect to Wuhan \cite{fays2}. Precisely, 
the fractional-order model proposed in \cite{fays2} describes well the spread of 
COVID-19 in Spain with order $\alpha =0.85$, in Portugal with order $\alpha = 0.75$, 
and in Wuhan with order $\alpha=1$. 

Fractional differential equations can be seen as a sub-field of fractional calculus \cite{samko1993fractional}, 
that is, the mathematical theory that deals with generalization of integrals and derivatives 
to real or complex order. Furthermore, systems modeled with the help of fractional calculus 
are non-linear and might display a much more richer dynamical behavior due to properties of order differentiation.  

The paper is organized as follows. We start by a preliminary section (Section~\ref{sec:prelim}), 
in which we recall definitions and necessary results of the literature. 
Then, in Section~\ref{sec:model}, we present the fractional-order model that we will be studying. 
Next, in Section~\ref{sec:exist:unic}, we prove the existence and uniqueness of a positive solution. 
The global stability analysis of the disease free equilibrium point is investigated 
in Section~\ref{sec:stability}. In Section~\ref{sec:numeric}, we present some numerical simulations 
of the model. We end with Section~\ref{sec:conc}, with our conclusions.


\section{Preliminaries on Fractional Calculus}
\label{sec:prelim}

In this section, we begin by presenting the Caputo definition of fractional
derivative and then recall some basic  properties useful to study the 
fractional-order model. The reader interested to learn about fractional calculus is referred 
to \cite{MR2218073,MR1658022}.

The Caputo fractional derivative of order $\alpha\in(0,1)$ 
of a function ${x:[0,+\infty)\rightarrow\mathbb R}$ is given by
\[
{}^{C}D^{\alpha}x(t)=\frac{1}{\Gamma(1-\alpha)} 
\int_{0}^{t} (t-s)^{-\alpha} x'(s)ds,
\]
where $\Gamma(1-\alpha)=\displaystyle \int_0^{\infty}t^{-\alpha}\exp(-t)dt$ 
is the Euler Gamma function. 

Note that the value of the Caputo fractional derivative of the function $x$ at point $t$ involves all
the values of $x'(s)$ for $s \in [0,t]$ and, hence, it incorporates the history of $x$. We  also see that
${}^{\texttt C}D^{\alpha}x(t)$ tends to $x'(t)$ as $\alpha \to 1$. The next result is crucial in the
study of an initial fractional-order value problem, useful in the proof of our Theorem~\ref{theo:uniq}. 

Let $f: \mathbb{R}^n \rightarrow \mathbb{R}^n$ be a vector function with $n > 1$ 
and consider the following fractional-order initial value problem:
\begin{equation}
\label{syst:exis}
\begin{cases}
{}^{C}D^{\alpha}X(t)=f(X),\\[3mm]
X(0)=X_0, \, \, X_0\in  \mathbb{R}^n.
\end{cases}
\end{equation}

\begin{Lemma}[See \cite{MR4232864}]
Assume that the vector function $f$ satisfies the following conditions:
\begin{enumerate}
\item $f(X)$ and $\displaystyle{\frac{\partial f(X)}{\partial X}}$ 
are continuous for all $X\in\mathbb{R}^n$;

\item $\Vert f(X) \Vert \leq \omega +\lambda \Vert X\Vert$
for all $X\in \mathbb{R}^n$, where $\omega$ and $\lambda$ 
are two positive constants.
\end{enumerate}

Then, system \eqref{syst:exis} has a unique solution.
\end{Lemma}

The following generalized mean value theorem and its consequences 
are also needed in the proof of Theorem~\ref{theo:uniq}.

\begin{Lemma}[Generalized Mean Value Theorem \cite{MR4232864}]
\label{mean:v}
Suppose that the functions $x(t)$ and ${}^{C}D^{\alpha}x(t)$ 
are both continuous on $[0, b]$. Then, 
\[
x(t)= x(0)+ \frac{1}{\Gamma(\alpha)}{}^{C}D^{\alpha}x(\eta)t^{\alpha}, 
\quad 0<\eta <t, \, \, \forall t\in [0, b].
\]
\end{Lemma}

Thus, as consequences of Lemma~\ref{mean:v}, 
we have that if ${}^{C}D^{\alpha}x(t)>0$ for all $t\in [0, b]$, 
then the function $x$ is strictly increasing, and if ${}^{C}D^{\alpha}x(t)<0$ 
for all  $t\in [0, b]$, then the function $x$ is strictly decreasing.


\section{The Considered Fractional-Order COVID-19 Model}
\label{sec:model}

In this section, we consider a fractional-order COVID-19 model earlier proposed 
by Nda\"{\i}rou et al. \cite{fays2}. We shall assume the total population 
$N$ is constant, along the period under study, and made up with eight sub-population of
dynamics transition, as different stages of transmission of the virus to individuals 
grouped into compartmental classes
$$
S(t)+E(t)+I(t)+P(t)+A(t)+H(t)+R(t)+F(t)=N 
$$
for all $t$, where $S(t)$ denotes the susceptible individuals at time $t$,
$E(t)$ the exposed individuals, $I(t)$ the symptomatic and infectious individuals, 
$P(t)$ the super-spreaders individuals, $A(t)$ the infectious but asymptomatic individuals, 
$H(t)$ the hospitalized individuals, $R(t)$ the recovery individuals, and 
$F(t)$ the dead individuals or fatality class. The Caputo fractional-order system that describes 
the dynamics transmission is given by
\begin{equation}
\label{model}
\begin{cases}
\displaystyle{{}^{\textsc c}D^{\alpha}S(t) 
= -\beta\frac{I(t)}{N}S(t)-l\beta \frac{H(t)}{N}S(t)
-\beta^{'}\frac{P(t)}{N}S(t)},\\[3mm]
\displaystyle{{}^{\textsc c}D^{\alpha}E(t)
= \beta\frac{I(t)}{N}S(t)+l\beta \frac{H(t)}{N}S(t)
+ \beta^{'}\frac{P(t)}{N}S(t) -\kappa E(t)}, \\[3mm]
\displaystyle{{}^{\textsc c}D^{\alpha}I(t)
= \kappa \rho_1 E(t) - (\gamma_a + \gamma_i)I(t)-\delta_i I(t)}, \\[3mm]
\displaystyle{{}^{\textsc c}D^{\alpha}P(t)
= \kappa \rho_2 E(t)- (\gamma_a + \gamma_i)P(t)
-\delta_p P(t)}, \\[3mm]
\displaystyle{{}^{\textsc c}D^{\alpha}A(t)
= \kappa (1-\rho_1 - \rho_2)E(t)},\\[3mm]
\displaystyle{{}^{\textsc c}D^{\alpha}H(t)
= \gamma_a (I(t) + P(t)) - \gamma_r H(t) 
- \delta_h H(t)}, \\[3mm]
\displaystyle{{}^{\textsc c}D^{\alpha}R(t)
= \gamma_i (I(t) + P(t))+ \gamma_r H(t)},\\[3mm]
\displaystyle{{}^{\textsc c}D^{\alpha}F(t)
= \delta_i I(t) + \delta_p P(t) + \delta_h H(t)}.
\end{cases}
\end{equation}

The expression $\beta\frac{I}{N}S+l\beta \frac{H}{N}S+ \beta^{'}\frac{P}{N}S$ 
represents the force of infection of the virus, that is, the transmission term 
or the effective contact between susceptible individuals ($S$) and infectious 
symptomatic individuals ($I$), super-spreaders individuals ($P$), and hospitalized 
ones ($H$). Here, $\beta$ quantifies the human-to-human transmission 
coefficient per unit of time (days) per person, $\beta^{'}$ quantifies a high 
transmission coefficient due to super-spreaders, $l$ quantifies the relative 
transmissibility of hospitalized patients. Next, we give a description of the 
rest of parameters appearing in the model system \eqref{model}:
\begin{itemize}
\item $\kappa$ is the rate at which an individual leaves the exposed 
class by becoming infectious (symptomatic, super-spreaders or asymptomatic);
\item $\rho_1$ is the proportion of progression from exposed class $E$ to symptomatic infectious class $I$;
\item $\rho_2$ is a relative very low rate at which exposed individuals become super-spreaders;
\item $1-\rho_1-\rho_2$ is the progression from exposed to asymptomatic class;
\item $\gamma_a$ is the average rate at which symptomatic  
and super-spreaders individuals become hospitalized;
\item $\gamma_i$ is the recovery rate without being~hospitalized;
\item $\gamma_r$ is the recovery rate of hospitalized patients;
\item $\delta_i$ denotes the disease induced death rates due to infected individuals; 
\item $\delta_p$ denotes the disease induced death rates due to super-spreaders individuals;
\item $\delta_h$ denotes the disease induced death rates due to hospitalized individuals.
\end{itemize}

The model was first proposed in \cite{fays2} for the purpose
of fitting real data from Galicia, Spain, and Portugal,
but without any mathematical analysis. Here we show
that the model is mathematically well posed and it
has a unique equilibrium point, which is globally
asymptotically stable.


\section{Existence and Uniqueness of Positive Solution}
\label{sec:exist:unic}

First of all, let us rewrite system \eqref{model} in a compact form. 
In doing so, denote 
$$
\mathbb{R}^8_{+}= \{X\in \mathbb{R}^8: X\geqslant 0\}
$$ 
and let $X(t)= \left(S(t), E(t), I(t), P(t), A(t), H(t), R(t), F(t) \right)^{T}$. 
Then, the system \eqref{model} can be rewritten as follows:
$$
{}^{\textsc c}D^{\alpha}X(t)= F\left(X(t) \right), 
$$
where
\begin{equation}
F(X)=\begin{pmatrix}
-\beta\frac{I}{N}S-l\beta \frac{H}{N}S-\beta^{'}\frac{P}{N}S\\
\beta\frac{I}{N}S+l\beta \frac{H}{N}S+ \beta^{'}\frac{P}{N}S -\kappa E\\
\kappa \rho_1 E - (\gamma_a + \gamma_i)I-\delta_i I\\
\kappa \rho_2 E- (\gamma_a + \gamma_i)P-\delta_p P\\
 \kappa (1-\rho_1 - \rho_2)E\\
 \gamma_a (I + P) - \gamma_r H - \delta_h H\\
 \gamma_i (I + P)+ \gamma_r H\\
 \delta_i I + \delta_p P + \delta_h H
\end{pmatrix}.
\end{equation}

For practical applications reasons, we consider non-negative initial conditions:
\begin{equation}
\label{cond}
S(0)\geqslant 0,\, \, E(0)\geqslant 0, \, \,  I(0)\geqslant 0, \, \, 
P(0)\geqslant 0,\, \, A(0)\geqslant 0, \, \, H(0)\geqslant 0, \, \, 
R(0)\geqslant 0, \, \, F(0)\geqslant 0.
\end{equation}

In addition, set
\begin{gather*}
A_1= \begin{pmatrix}
0&0&-\beta & -\beta^{'}&0&-l\beta &0&0\\
0&0&\beta & \beta^{'}&0&l\beta &0&0\\
0&0&0&0&0&0&0&0\\
0&0&0&0&0&0&0&0\\
0&0&0&0&0&0&0&0\\
0&0&0&0&0&0&0&0\\
0&0&0&0&0&0&0&0\\
0&0&0&0&0&0&0&0\\
\end{pmatrix}, \\
A_2=\begin{pmatrix}
0&0&0&0&0&0&0&0\\
0&-\kappa&0&0&0&0&0&0\\
0&\kappa \rho_1&-\varpi_i&0&0&0&0&0\\
0&\kappa \rho_2&0&\varpi_p&0&0&0&0\\
0&\varpi_e&0&0&0&0&0&0\\
0&0&\gamma_a&\gamma_a&0&\varpi_h&0&0\\
0&0&\gamma_i&\gamma_i&0&\gamma_r&0&0\\
0&0&0&\delta_i&\delta_p&0&\delta_h&0\\
\end{pmatrix},
\end{gather*}
where $\varpi_e = \kappa (1-\rho_1-\rho_2)$; 
$\varpi_i = \gamma_a + \gamma_i + \delta_i$; 
$\varpi_p= \gamma_a + \gamma_i + \delta_p$; 
and $\varpi_h= \gamma_r + \delta_h$.
Thus, we can rewrite the vector function $F$ as 
$$
F(X)= \frac{S}{N}A_1X + A_2X.
$$

We are in conditions to state and prove our first result.

\begin{Theorem}[existence and uniqueness of a non-negative solution]
\label{theo:uniq}
There is a unique solution for the initial value problem given by 
\eqref{model}--\eqref{cond} and the solution belongs to 
$$
\Omega = \{(S, E, I, P, A, H, R, F) 
\in \mathbb{R}^8_{+}: S+E+I+P+A+H+R+F\leqslant N \}. 
$$
\end{Theorem}

\begin{proof}
The existence and uniqueness of solution are obtained by application of Theorem~3.1 
and Remark~3.2 of \cite{lin}. For this purpose, it is easy to check that the vector 
function is a polynomial, thus continuous  and has continuous derivative in $\Omega$. 
Furthermore, satisfying
\[
\Vert F(X)\Vert \leqslant \Vert A_1 X \Vert + \Vert A_2 X\Vert 
= \left(\Vert A_1 \Vert + \Vert A_2\Vert\right)\Vert X \Vert < \epsilon 
+ \left(\Vert A_1 \Vert + \Vert A_2\Vert\right)\Vert X \Vert
\]
for any positive constant $\epsilon$. The proof of non-negativity of solution 
follows the same spirit as in \cite{ricardo}. By summing up all the 8 equations 
of system \eqref{model}, we obtain
\[
{}^{\textsc c}D^{\alpha} \left(S(t)+ E(t)+ I(t)
+ P(t)+ A(t)+ H(t)+ R(t)+ F(t) \right)=0.
\]

Thus,
\begin{align*}
0&\leqslant S(t)+ E(t)+ I(t)+ P(t)+ A(t)+ H(t)+ R(t)+ F(t)\\
&\leqslant S(0)+ E(0)+ I(0)+ P(0)+ A(0)+ H(0)+ R(0)+ F(0) =N,
\end{align*} 
which ends the proof.
\end{proof}


\section{Stability Analysis}
\label{sec:stability}

First of all, note that the model system \eqref{model} exhibits a unique 
steady state which is the disease free equilibrium point (DFE) obtained by 
setting the right hand side of \eqref{model} equal to zero. Precisely, we have 
$$
DFE = (N,0,0,0,0,0,0,0).
$$

Next, recall that the basic reproduction number for this fractional-order model system 
is the same as for the classical model investigated in \cite{fays1,MyID:460} 
using the next-generation matrix approach \cite{van:den:Driessche:2002, diekmann}, 
being given by
\begin{equation}
\label{R0}
R_0= \frac{\beta \rho_1(\gamma_a l + \varpi_h)}{\varpi_i  \varpi_h} 
+ \frac{(\beta \gamma_a l + \beta^{'}\varpi_h) \rho_2}{\varpi_p \varpi_h}.
\end{equation}

This can be rewritten in the following manner:
\begin{equation}
\label{R0:rewrite}
R_0= \frac{\beta\rho_1 \varpi_h\varpi_p + \beta \rho_1\gamma_al\varpi_p 
+ \beta^{'}\rho_2\varpi_h\varpi_i + \beta \rho_2\gamma_a 
l\varpi_i}{\varpi_i \varpi_p \varpi_h},
\end{equation}
which is useful in the below proof of global stability. 

\begin{Theorem}[global stability of the DFE]
Let $\alpha \in (0, 1)$. The disease free equilibrium (DFE)
of system \eqref{model} is globally asymptotically stable whenever $R_0 < 1$. 
\end{Theorem}

\begin{proof}
Consider the following Lyapunov function:
\[
V(t)= a_0E(t)+a_1I(t)+a_2P(t)+a_3H(t),
\]
where $a_0$, $a_1$, $a_2$, and $a_3$ are positive constants to be determined.
Because the fractional operator ${}^{\textsc c}D^{\alpha}$ is linear, we have that 
\[
{}^{\textsc c}D^{\alpha} V(t)= a_0 {}^{\textsc c}D^{\alpha}E(t) 
+ a_1 {}^{\textsc c}D^{\alpha}I(t) + a_2 {}^{\textsc c}D^{\alpha}P(t) 
+ a_3{}^{\textsc c}D^{\alpha}H(t),
\]
and from \eqref{model} it follows that
\begin{align*}
{}^{\textsc c}D^{\alpha}V(t)&= a_0\left(\beta\frac{I}{N}S
+l\beta \frac{H}{N}S+ \beta^{'}\frac{P}{N}S -\kappa E \right) 
+ a_1\left(\kappa \rho_1 E - (\gamma_a + \gamma_i)I-\delta_i I \right) \\
&\quad +a_2\left(\kappa \rho_2 E- (\gamma_a + \gamma_i)P-\delta_p P \right) 
+ a_3\left( \gamma_a (I + P) - \gamma_r H - \delta_h H\right).
\end{align*}

Next, as $S\leqslant N$, we have 
\begin{align*}
{}^{\textsc c}D^{\alpha}V(t)&\leqslant a_0\left(\beta I+l\beta H+ \beta^{'}P -\kappa E \right) 
+ a_1\left(\kappa \rho_1 E - \varpi_i I \right) \\
&\quad +a_2\left(\kappa \rho_2 E- \varpi_p P \right) 
+ a_3\left( \gamma_a (I + P) - \varpi_h H\right),
\end{align*}
where $\varpi_i = \gamma_a + \gamma_i + \delta_i$;
$\varpi_p= \gamma_a + \gamma_i + \delta_p$; 
and $\varpi_h= \gamma_r + \delta_h$.
Rearranging and reducing leads to 
\begin{align*}
{}^{\textsc c}D^{\alpha}V(t)&\leqslant  (a_0\beta 
+ a_3\gamma_a -a_1\varpi_i)I + (a_0\beta l -a_3\varpi_h)H \\
&\quad + (a_0\beta^{'} + a_3\gamma -a_2\varpi_p)P + \kappa(a_1\rho_1 + a_2\rho_2-a_0)E.
\end{align*}

Now, we choose,
\begin{gather*}
a_0= \varpi_i \varpi_p \varpi_h; 
\quad a_1 = \left(\beta + \frac{\beta \gamma_a l}{\varpi_h} \right)\varpi_h \varpi_p; 
\quad a_2= \left(\beta^{'} + \frac{\beta \gamma_a l}{\varpi_h} \right)\varpi_i \varpi_h; 
\quad a_3= \beta l \varpi_i \varpi_p,
\end{gather*}
so that function $V$ is defined, continuous, and positive definite for all 
$E(t)>0$, $I(t)>0$, $P(t)>0$, and $H(t)>0$. As a consequence, we obtain that 
\begin{gather*}
a_0\beta + a_3\gamma_a -a_1\varpi_i = 0; 
\quad a_0\beta l -a_3\varpi_h= 0; 
\quad a_0\beta^{'} + a_3\gamma -a_2\varpi_p=0,
\end{gather*} 
and 
\begin{align*}
a_1\rho_1 + a_2\rho_2-a_0 
&= \beta \rho_1\varpi_h\varpi_p + \beta \rho_1\gamma_a l\varpi_p + \beta^{'}\rho_2\varpi_h\varpi_i 
+ \beta \rho_2\gamma_a l\varpi_i-\varpi_i \varpi_p \varpi_h\\
&= \varpi_i \varpi_p \varpi_h \left(\frac{\beta\rho_1 \varpi_h\varpi_p + \beta \rho_1\gamma_al\varpi_p 
+ \beta^{'}\rho_2\varpi_h\varpi_i + \beta \rho_2\gamma_a l\varpi_i}{\varpi_i \varpi_p \varpi_h} 
-1 \right).
\end{align*}

Note that from \eqref{R0:rewrite} we have that  
\[
R_0= \frac{\beta\rho_1 \varpi_h\varpi_p + \beta \rho_1\gamma_al\varpi_p 
+ \beta^{'}\rho_2\varpi_h\varpi_i + \beta 
\rho_2\gamma_a l\varpi_i}{\varpi_i \varpi_p \varpi_h}, 
\]
which by substitution leads to
\[
{}^{\textsc c}D^{\alpha} V(t)\leqslant 
\kappa \varpi_i \varpi_p \varpi_h (R_0 -1)E.
\]

Finally, ${}^{\textsc c}D^{\alpha} V(t)\leqslant 0$ if $R_0 <1$. 
Furthermore, ${}^{\textsc c}D^{\alpha}V(t)=0$ if, and only if, $E=I=P=H=0$. 
Substituting $(E, I, P, H)=(0,0,0,0)$ in system \eqref{model}, 
leads to 
$$
S(t)=S(0), \quad A(t)=A(0), \quad R(t)=R(0), \quad F(t)=F(0).
$$ 

Thus, the largest compact invariant set containing the DFE is  
$$
\Gamma =\lbrace (S, E, I, P, A, H, R, F)
\in \mathbb{R}^8_{+}:{}^{\textsc c}D^{\alpha} V(t)=0\rbrace.
$$

However, from biological considerations, when $(E, I, P, H)=(0,0,0,0)$, 
meaning there is no disease infection in the population, we have the implication 
$A(0)=R(0)=F(0)=0$ and $S(0)=N$. Therefore, the largest compact invariant $\Gamma$ 
set is reduced to the singleton $\lbrace DFE \rbrace$. Hence, 
by LaSalle invariance principle \cite{MR3384337}, 
we conclude that the disease free equilibrium DFE 
is globally asymptotically stable.
\end{proof}

In previous paper \cite{fays2}, the case studies
of Galicia, Spain, and Portugal are investigated,
separately, with the purpose to fit the real data.
In the next section we use the same values of 
fractional order given in \cite{fays2} and we focus 
on a comparative study of the three mentioned cases
of Galicia, Spain, and Portugal, together with 
the Wuhan case studied in \cite{fays1,MyID:460}.


\section{Numerical Simulations}
\label{sec:numeric}

In this section, model analysis is carried out through numerical simulations 
in order to show a broad view of the time evolution of the infected populations. 
Mainly, we shall study the dynamical behavior of infected individuals ($I$), 
super-spreaders ($P$), hospitalized individuals ($H$), 
and the cumulative cases of infections ($I + P+ H$), obtained from the output 
of our fractional-order system \eqref{model}. We will focus on a comparative study 
by considering values of order of differentiation $\alpha=1$, 
$\alpha=0.85$, and $\alpha=0.75$, which describe the COVID-19 dynamics transmission 
of Wuhan, Spain, and Portugal, respectively \cite{fays2}. In addition, the role of 
the basic reproduction number through infectivity effect on the evolution curves 
will be conducted. The readers interested in seeing the real data are referred to~\cite{fays2}.


\subsection{Population Size, Initial Conditions, and Parameters}

Recall that the total population size under study reflects specificities 
on the spread of COVID-19 on each territories considered. Therefore, we consider
$N=$ 47,000,000/425, $N=$ 10,280,000/875, and $N=$ 11,000,000/250 for Spain, 
Portugal, and Wuhan, respectively. We remark that there was a typo in
\cite{fays2} for the value of $N$ in the case of Portugal.
The following initial conditions are considered: 
$$
S_0=47,000,000/425-11, E_0= 0, I_0= 10, P_0= 1, A_0= 0, H_0= 0, R_0= 0, F_0= 0,
$$
for Spain;
$$ 
S_0=10,280,000/875-5, E_0= 0, I_0= 4, P_0= 1, A_0= 0, H_0= 0, R_0= 0, F_0= 0,
$$
for Portugal; and 
$$
S_0=11,000,000/250-6, E_0= 0, I_0= 1, P_0= 5, A_0= 0, H_0= 0, R_0= 0, F_0= 0,
$$ 
for Wuhan. Moreover, the following values of parameters are borrowed 
from \cite{fays1,MyID:460}:
\begin{gather*}
\beta =2.55, \quad l =1.56, \quad \beta^{'}=7.65, 
\quad \kappa=0.25, \quad \rho_1=0.58, \quad \rho_2=0.001, \\ 
\gamma_a=0.94, \quad \gamma_i=0.27, \quad \gamma_r=0.5, 
\quad \delta_i=\delta_p=\delta_h= \frac{1}{23}.
\end{gather*} 


\subsection{Index of Memory's Influence}

\begin{figure}[H]
\widefigure
\subfloat[\footnotesize{Symptomatic/infectious individuals $I(t)$\\ 
for $\alpha \in \{ 0.75, 0.85, 1\}$}]{
\label{Inf:syn}\includegraphics[scale=0.6]{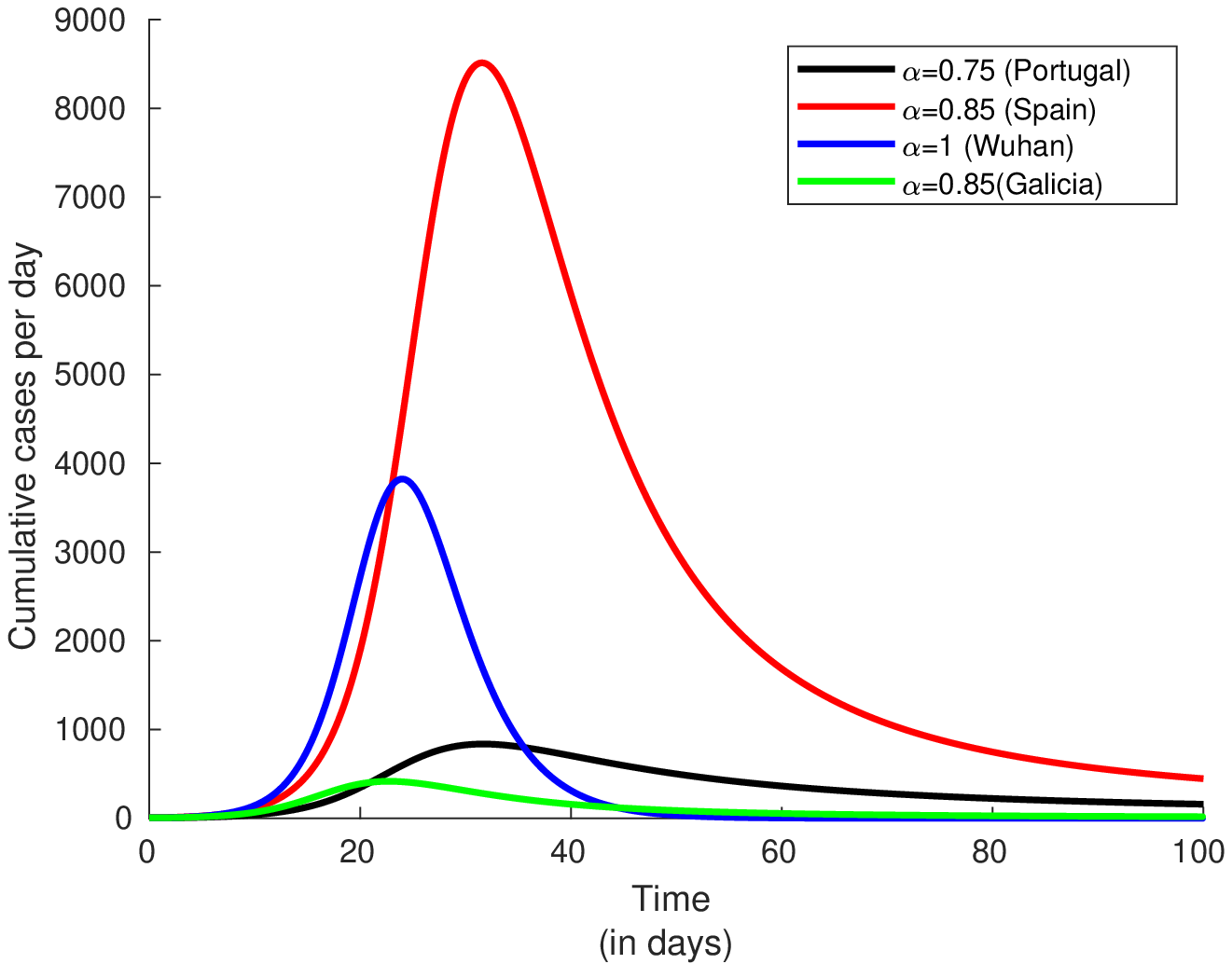}}
\subfloat[\footnotesize{Super-spreaders individuals $P(t)$\\ 
for $\alpha \in \{ 0.75, 0.85, 1\}$}]{
\label{Sspr:syn}\includegraphics[scale=0.6]{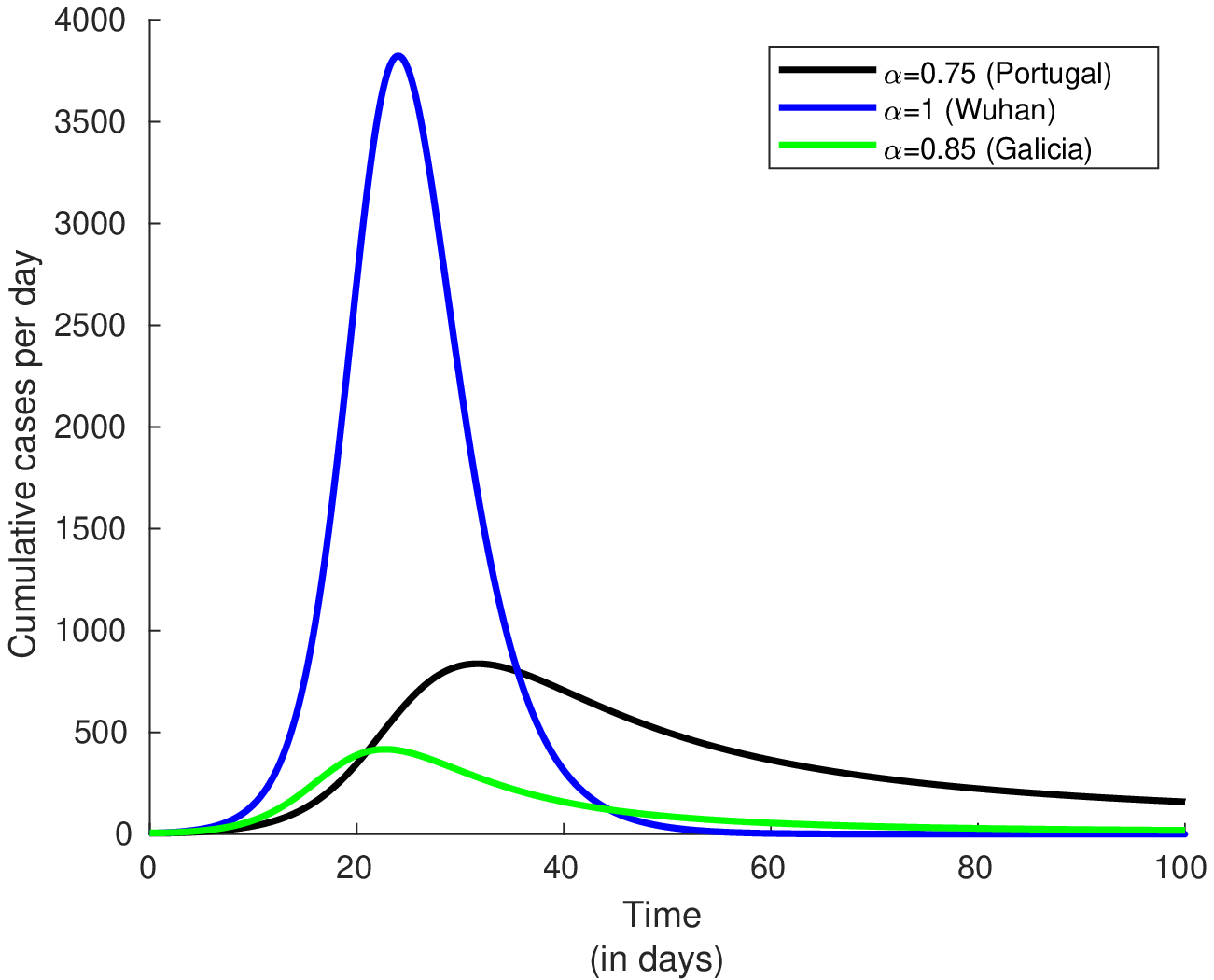}}\\
\subfloat[\footnotesize{Hospitalized individuals $H(t)$\\ 
for $\alpha \in \{ 0.75, 0.85, 1\}$}]{
\label{Hos:syn}\includegraphics[scale=0.6]{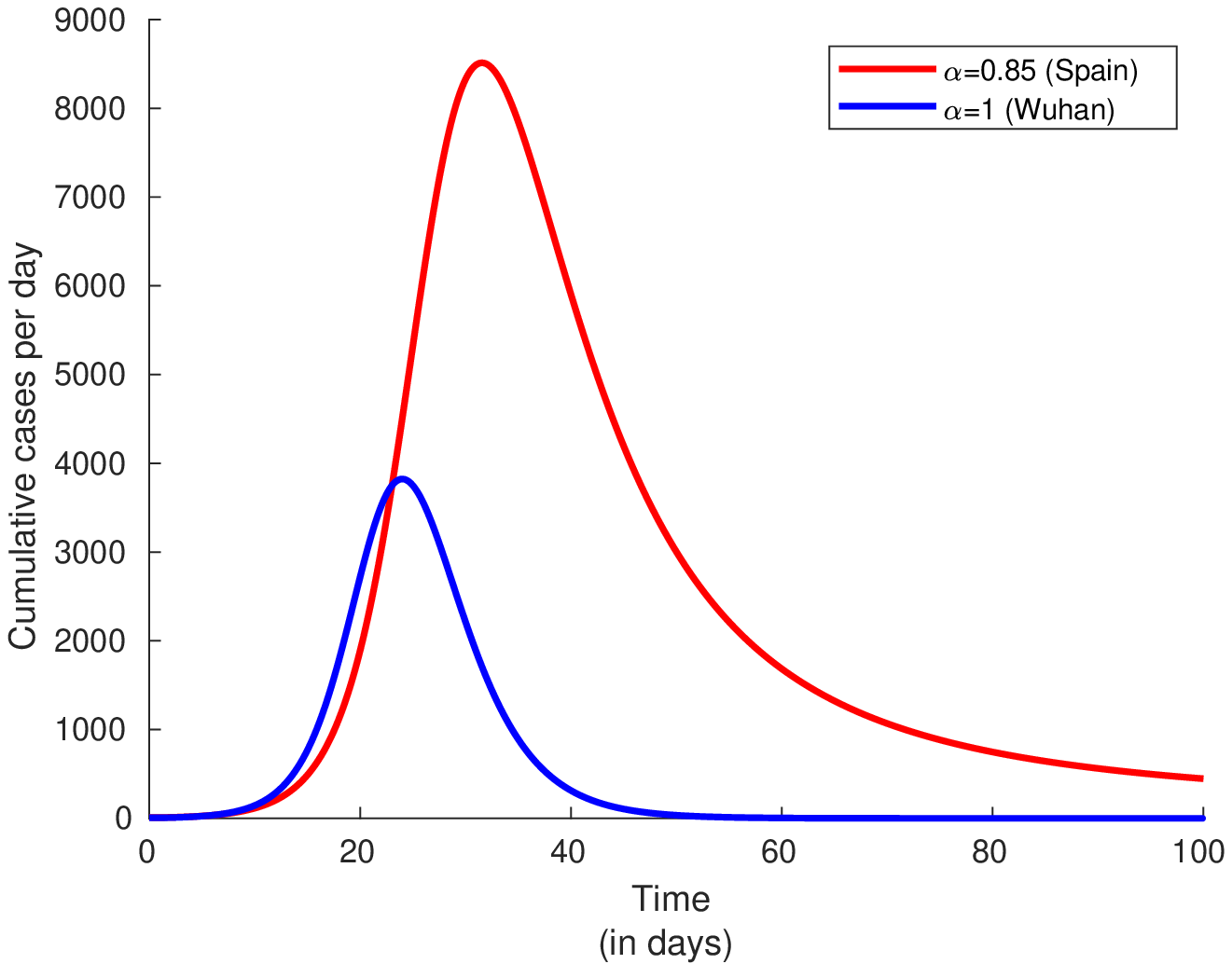}}
\subfloat[\footnotesize{$I(t)+ P(t) + H(t)$ for $\alpha \in \{ 0.75, 0.85, 1\}$}]{
\label{Cum:syn}\includegraphics[scale=0.6]{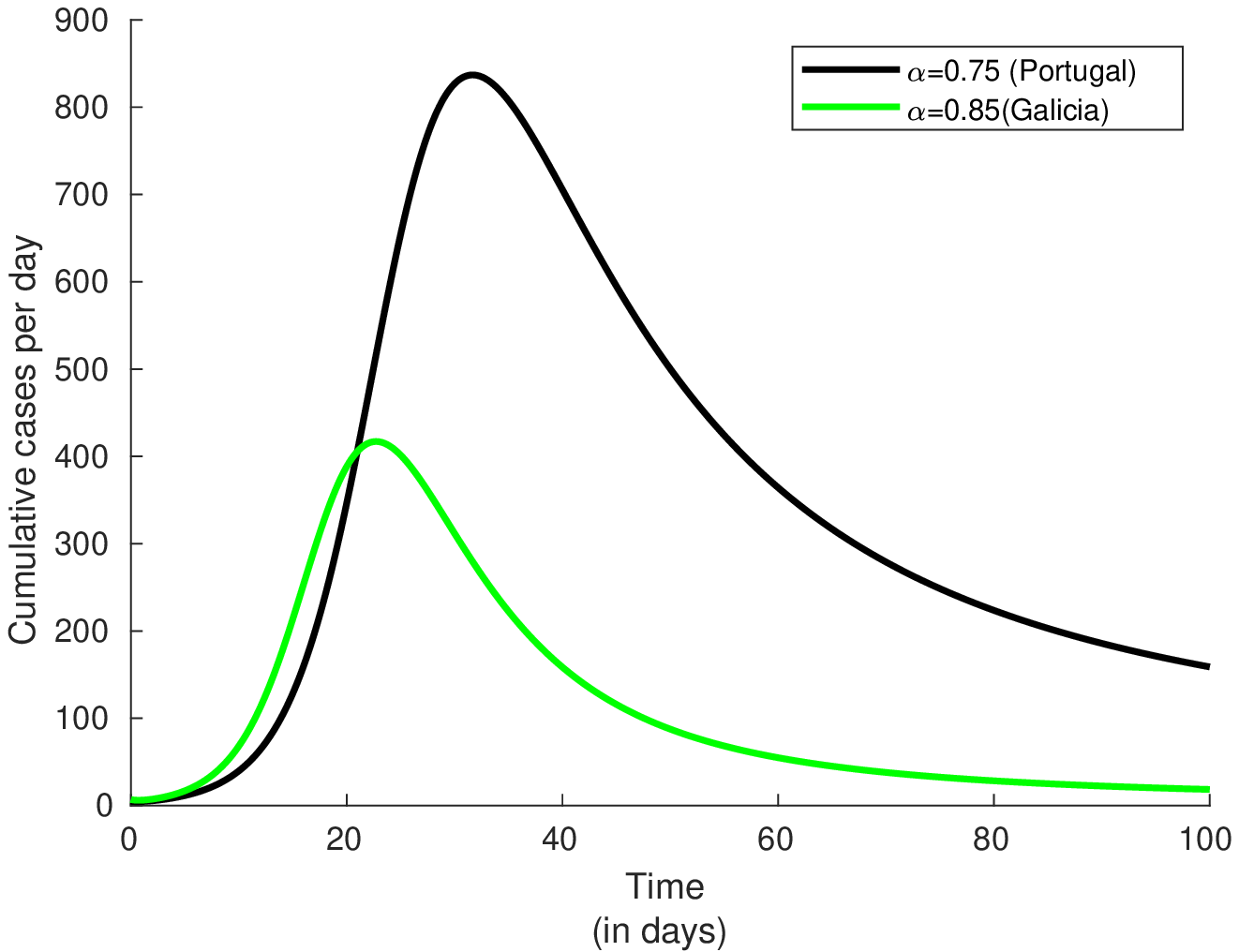}}
\caption{Evolution of infected populations by varying the order of differentiation $\alpha$:
the data of Portugal is well described with $\alpha = 0.75$; 
the data of Spain, and Galicia alone, with $\alpha = 0.85$; 
the data of Wuhan with $\alpha = 1$.}
\label{fig:evol:Inf}
\end{figure}

The stability of the cumulative cases of infections and the variation 
on the speed of convergence for different values of $\alpha$ are illustrated 
in Figure~\ref{fig:evol:Inf}. 


We observe that for a longer period of time, infected populations decrease 
and tend to zero. Further, the smaller the order of differentiation, 
the slower the convergence to the steady state. 

\begin{figure}[H]
\widefigure
\subfloat[\footnotesize{$I(t)$ for $\beta \in \{ 2.55, 1.55, 3.55\}$ but $\alpha=1$}]{
\label{Inf:wuhan}\includegraphics[scale=0.6]{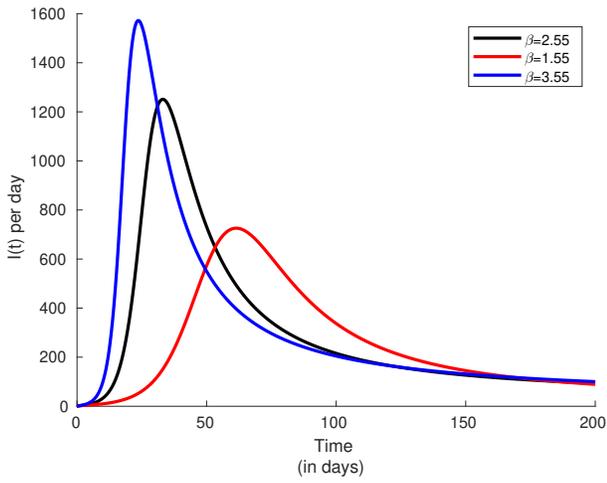}}
\subfloat[\footnotesize{$P(t)$ for $\beta \in \{ 2.55, 1.55, 3.55\}$ but $\alpha=1$}]{
\label{Sspr:wuhan}\includegraphics[scale=0.6]{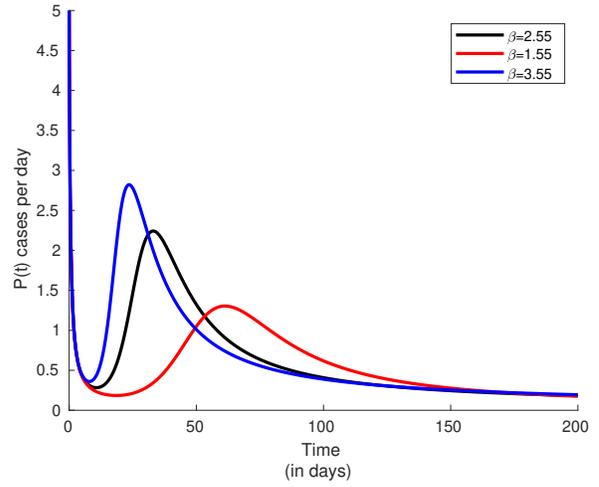}}\\
\subfloat[\footnotesize{$H(t)$ for $\beta \in \{ 2.55, 1.55, 3.55\}$ but $\alpha=1$}]{
\label{Hos:wuhan}\includegraphics[scale=0.6]{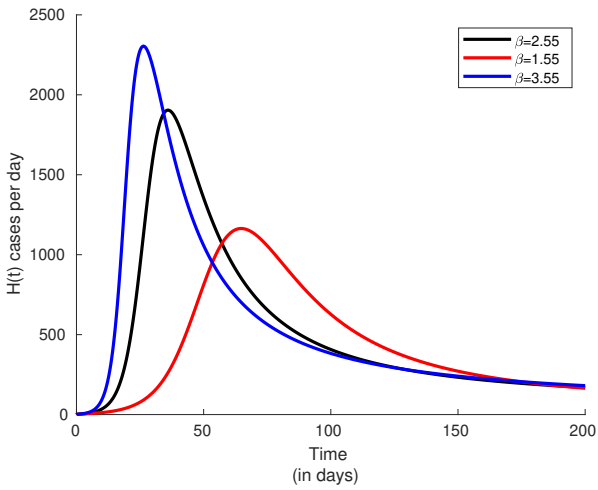}}
\subfloat[\footnotesize{$I(t)+ P(t) + H(t)$ for $\beta \in \{ 2.55, 1.55, 3.55\}$ but $\alpha=1$}]{
\label{Cum:wuhan}\includegraphics[scale=0.6]{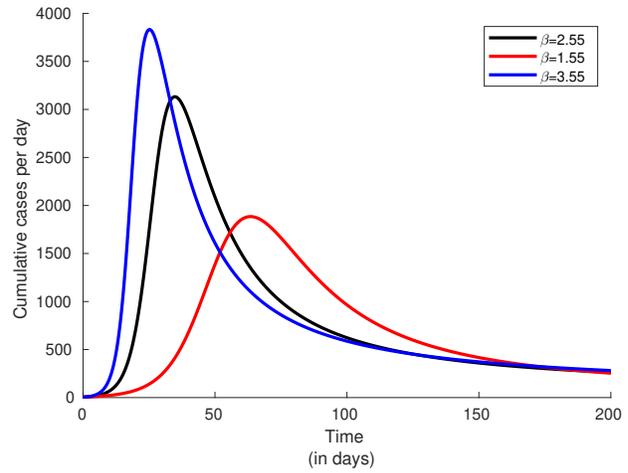}}
\caption{Evolution of infected populations ($I(t), P(t), H(t)$, and $ I(t)+P(t)+H(t)$) 
by varying the infectivity rate $\beta$ by $1.55$, $2.55$, and $3.55$, 
corresponding, respectively, to the basic reproduction number $2.662$, 
$4.375$, and $6.088$, while fixing index memory $\alpha=1$ (Wuhan).}
\label{fig:evol:wuhan}
\end{figure}
\begin{figure}[H]
\widefigure
\subfloat[\footnotesize{$I(t)$ for $\beta \in \{ 2.55, 1.55, 3.55\}$ but $\alpha=0.85$}]{
\label{Inf:spain}\includegraphics[scale=0.6]{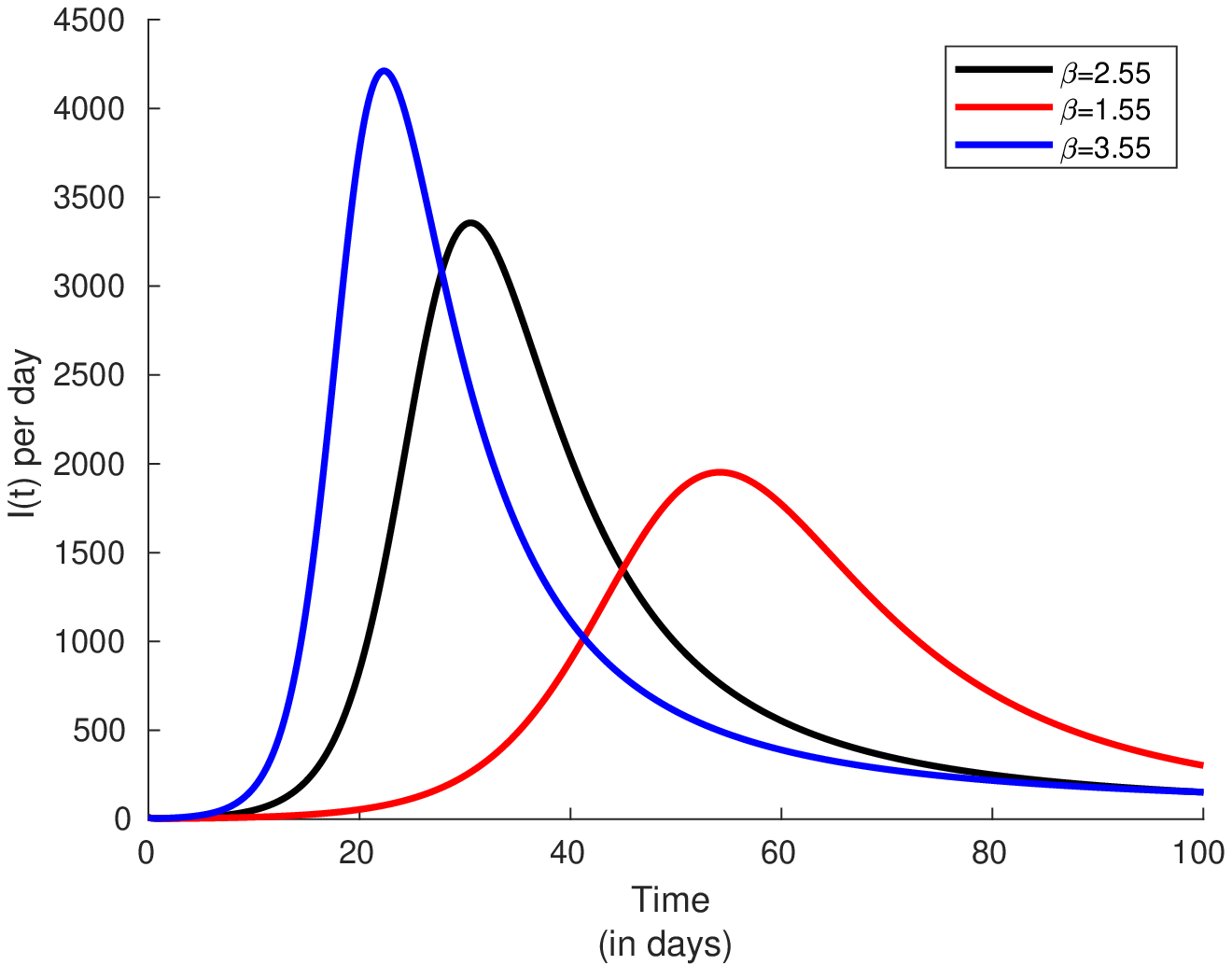}}
\subfloat[\footnotesize{$P(t)$ for $\beta \in \{ 2.55, 1.55, 3.55\}$ but $\alpha=0.85$}]{
\label{Sspr:spain}\includegraphics[scale=0.6]{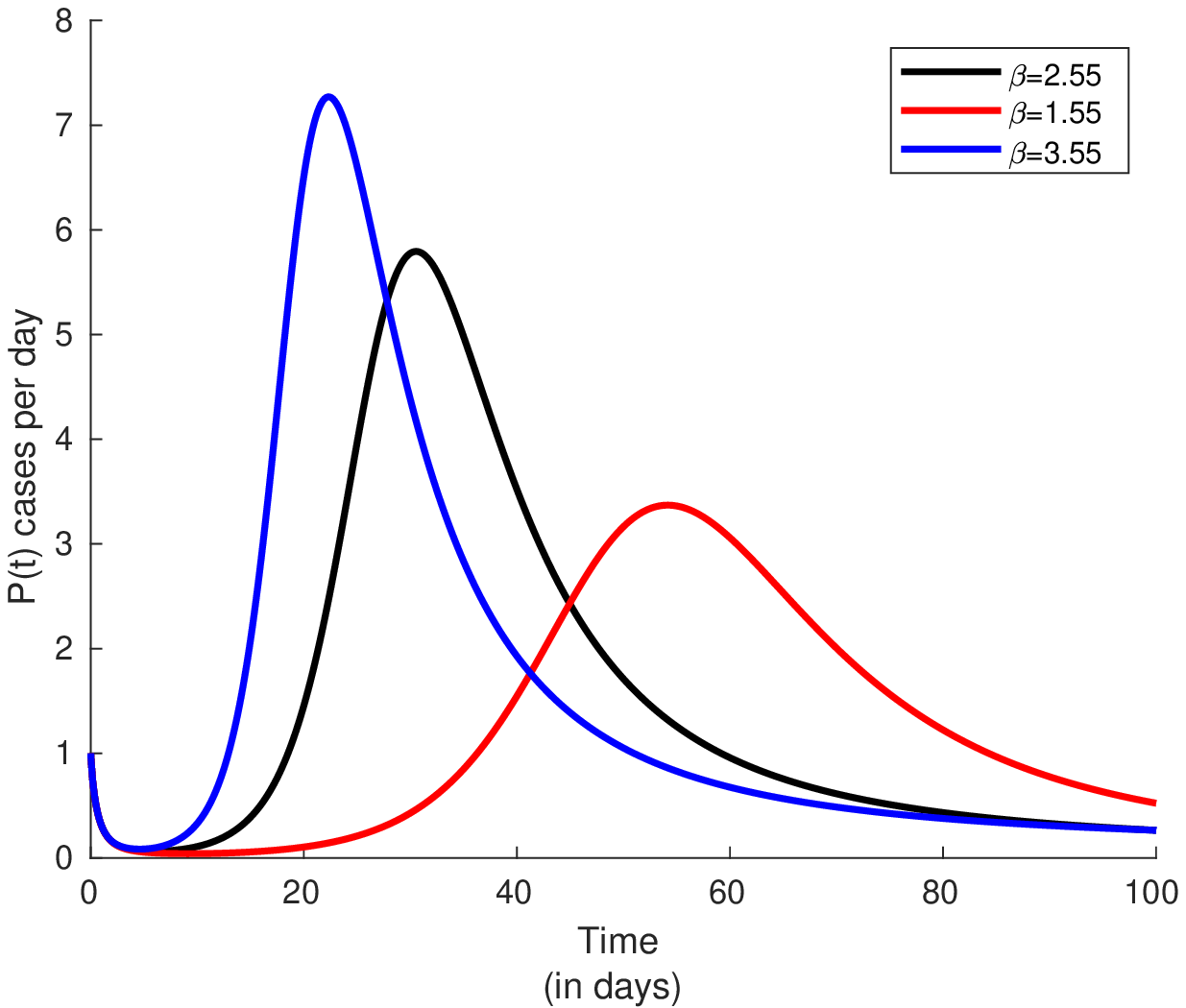}}\\
\subfloat[\footnotesize{$H(t)$ for $\beta \in \{ 2.55, 1.55, 3.55\}$ but $\alpha=0.85$}]{
\label{Hos:spain}\includegraphics[scale=0.6]{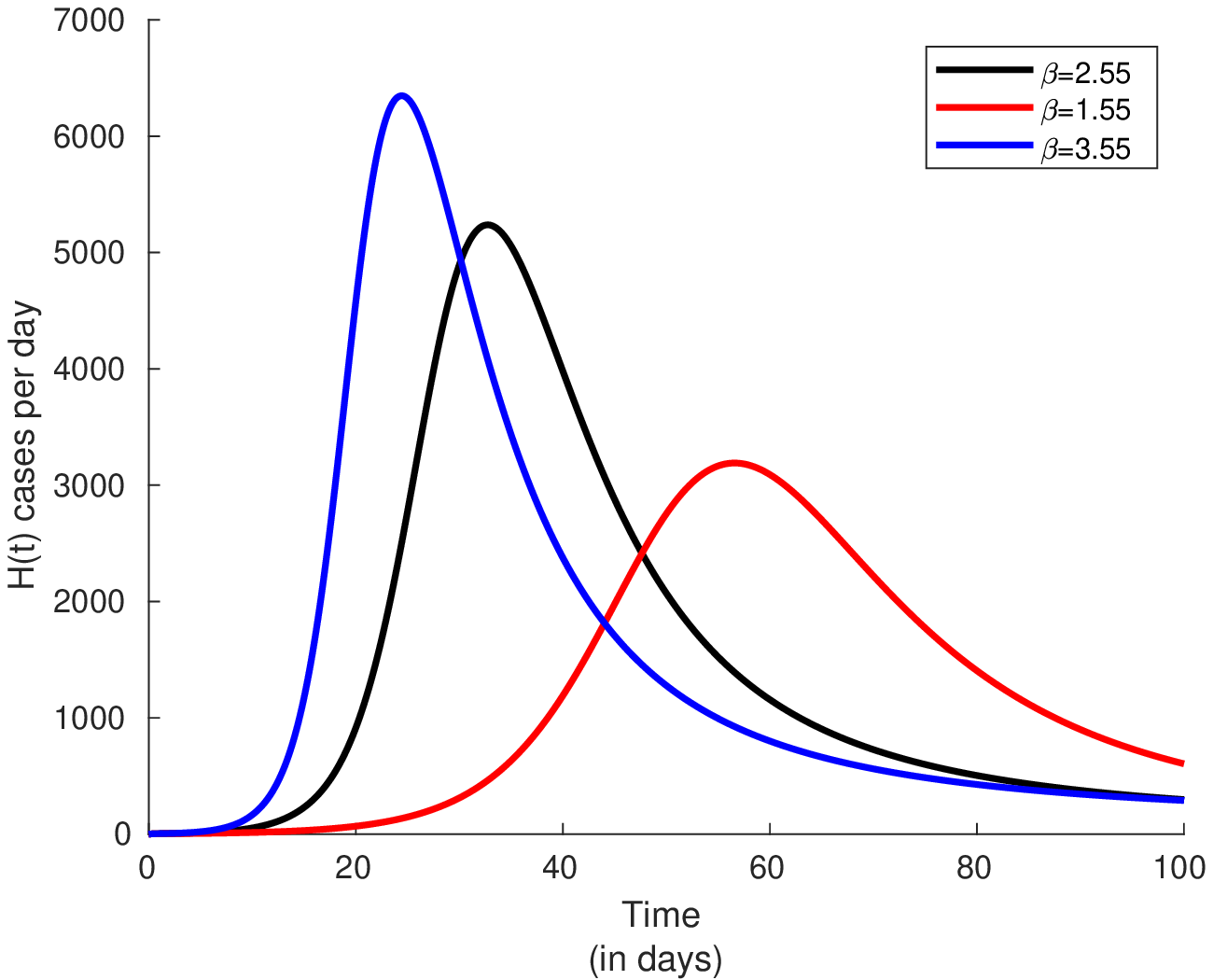}}
\subfloat[\footnotesize{$I(t)+ P(t) + H(t)$ for $\beta \in \{ 2.55, 1.55, 3.55\}$
\newline but $\alpha=0.85$}]{
\label{Cum:spain}\includegraphics[scale=0.6]{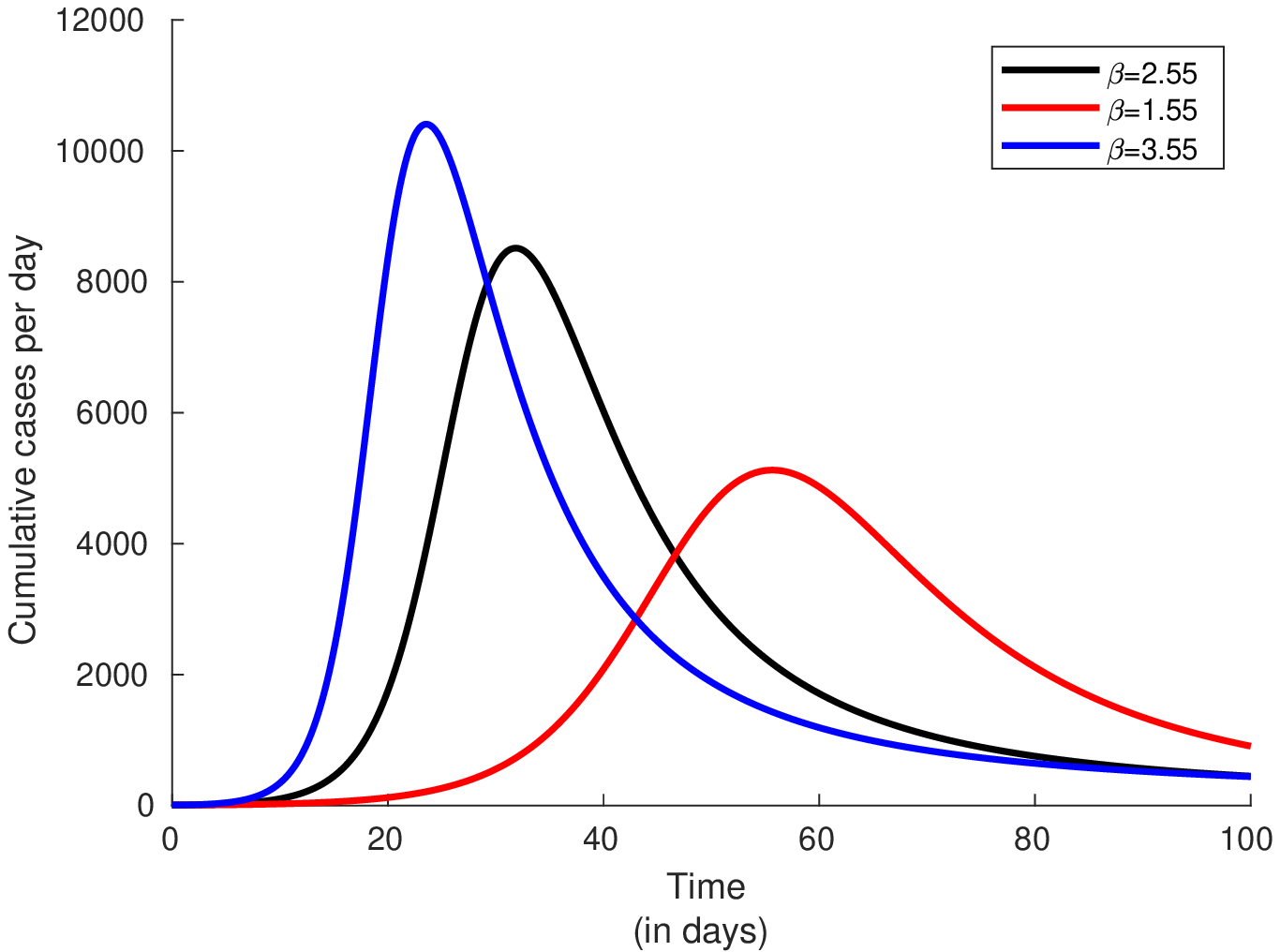}}
\caption{Evolution of infected populations ($I(t), P(t), H(t)$, and $ I(t)+P(t)+H(t)$) 
by varying the infectivity rate $\beta$ by $1.55$, $2.55$, and $3.55$, 
corresponding, respectively, to the basic reproduction number $2.662$, $4.375$, 
and $6.088$, while fixing index memory $\alpha=0.85$ (Spain).}
\label{fig:evol:spain}
\end{figure}
\begin{figure}[H]
\widefigure
\subfloat[\footnotesize{$I(t)$ for $\beta \in \{ 2.55, 1.55, 3.55\}$ but $\alpha=0.75$}]{
\label{Inf:portugal}\includegraphics[scale=0.6]{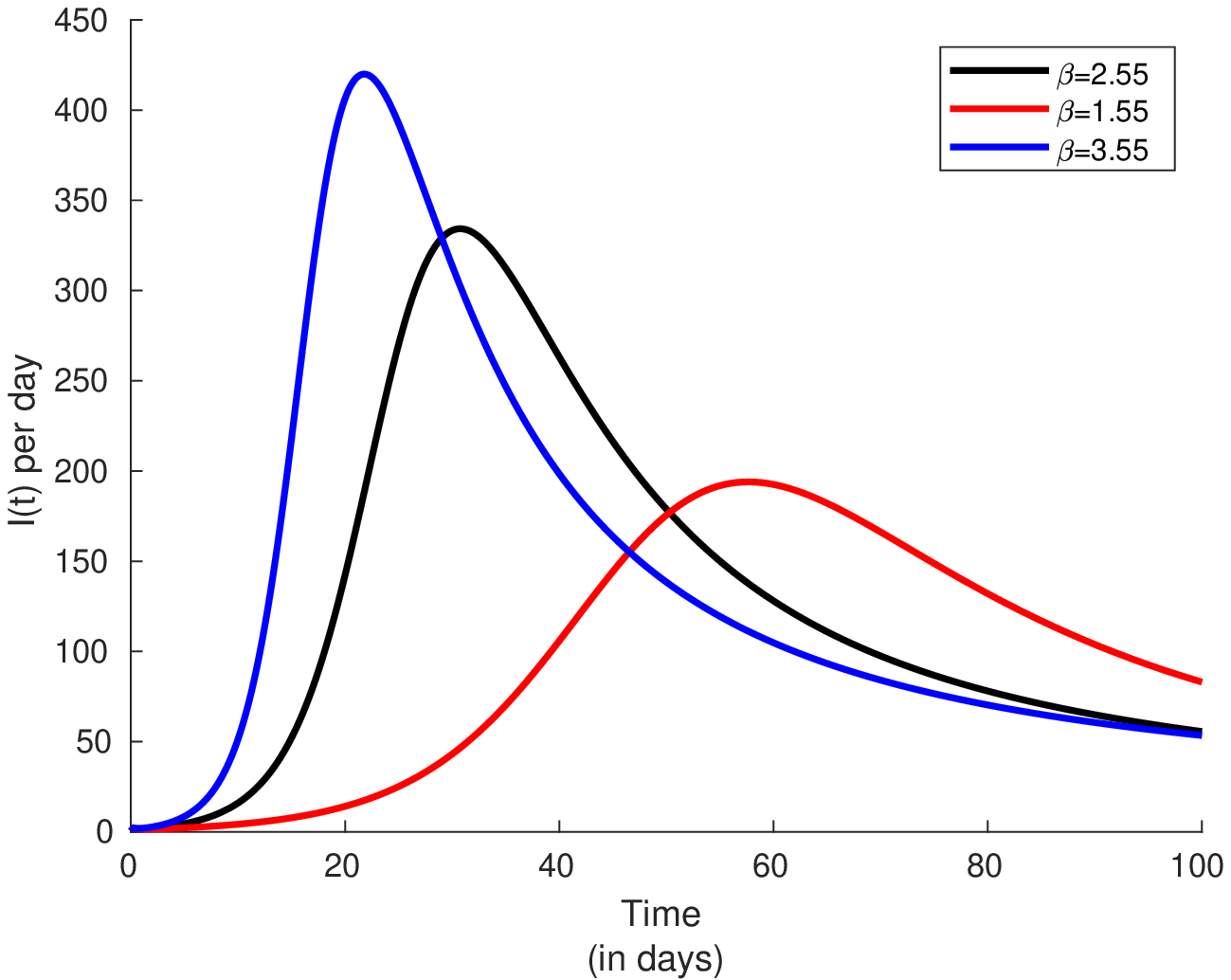}}
\subfloat[\footnotesize{$P(t)$ for $\beta \in \{ 2.55, 1.55, 3.55\}$ but $\alpha=0.75$}]{
\label{Sspr:portugal}\includegraphics[scale=0.6]{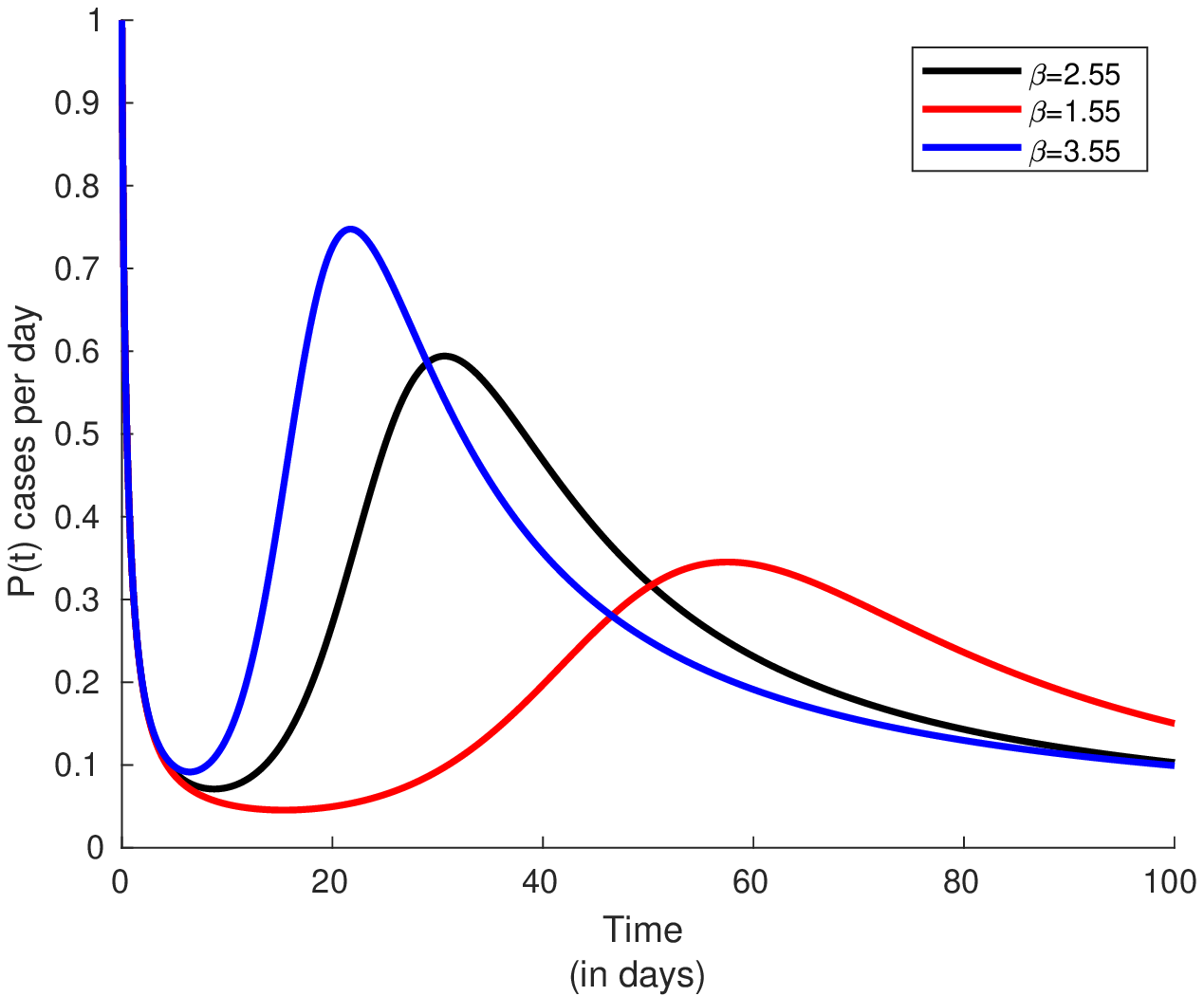}}\\
\subfloat[\footnotesize{$H(t)$ for $\beta \in \{ 2.55, 1.55, 3.55\}$ but $\alpha=0.75$}]{
\label{Hos:portugal}\includegraphics[scale=0.6]{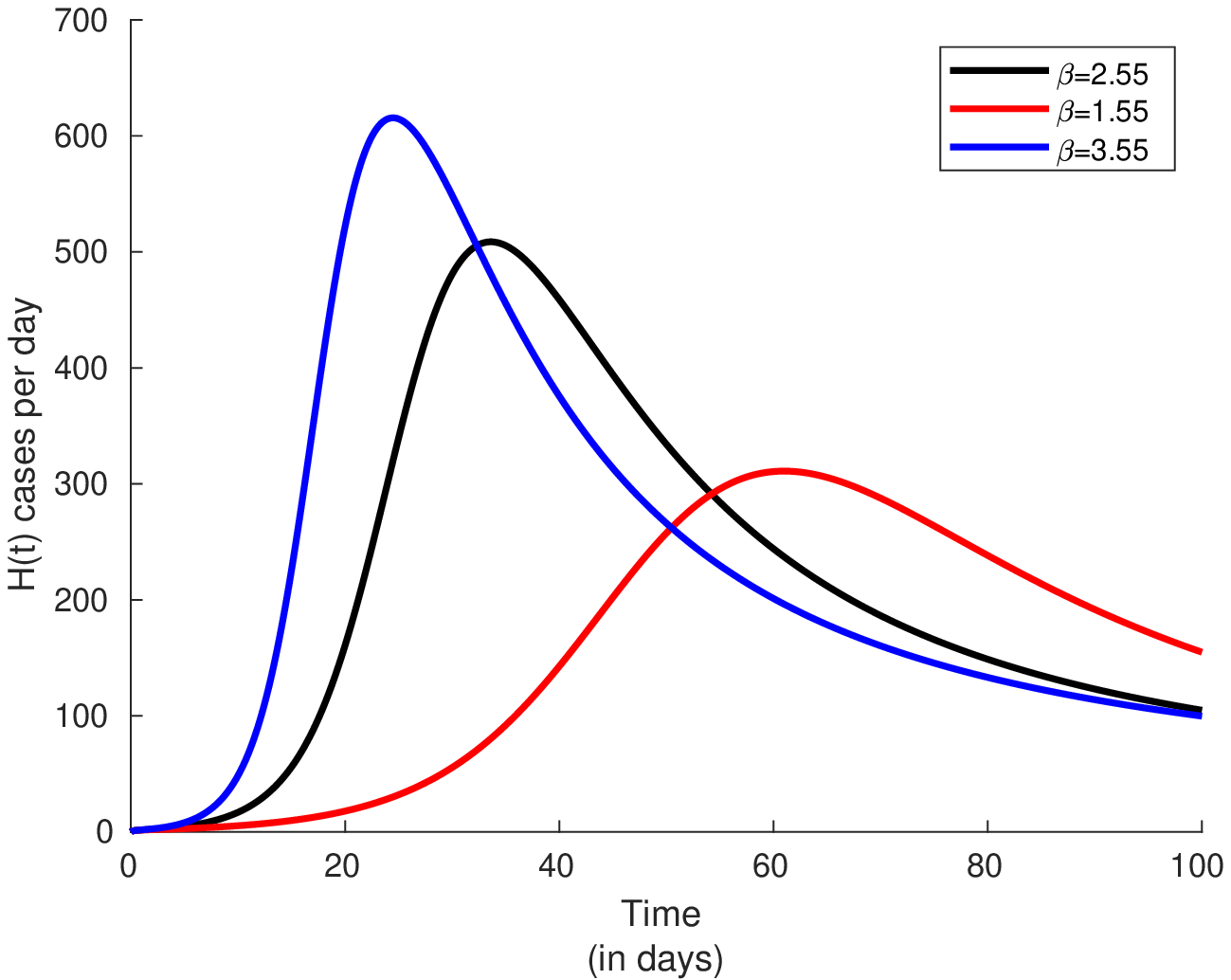}}
\subfloat[\footnotesize{$I(t)+ P(t) + H(t)$ for $\beta \in \{ 2.55, 1.55, 3.55\}$\newline 
but $\alpha=0.75$}]{
\label{Cum:portugal}\includegraphics[scale=0.6]{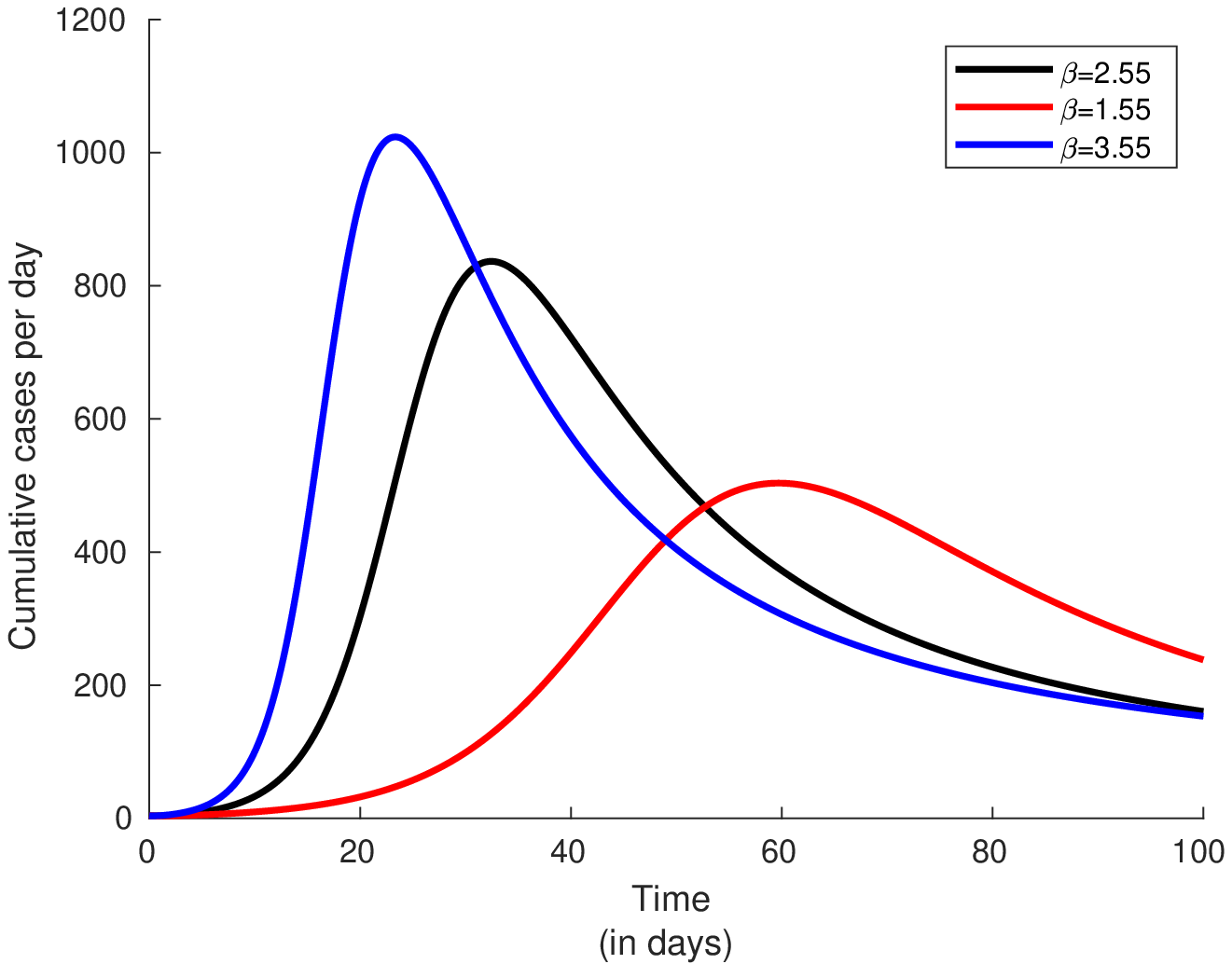}}
\caption{Evolution of infected populations ($I(t), P(t), H(t)$, and $ I(t)+P(t)+H(t)$) 
by varying the infectivity rate $\beta$ by $1.55$, $2.55$, and $3.55$, corresponding, 
respectively, to the basic reproduction number $2.662$, $4.375$, and $6.088$, 
while fixing index memory $\alpha=0.75$ (Portugal).}
\label{fig:evol:portugal}
\end{figure}

\clearpage
\end{paracol}
\nointerlineskip
\begin{paracol}{2}

\switchcolumn
\vspace{-6pt}


\subsection{Infectivity Rate and Effect on the Basic Reproduction Number}

To highlight the effect of the reproduction number $R_0$, from
Figures~\ref{fig:evol:wuhan}--\ref{fig:evol:portugal}, 
three different scenarios for the infectivity rate $\beta$ are considered 
with respect to each value of the index memory. Our results show that there 
is a significant decrease in the peak values of each infected categories 
of population when the reproduction number $R_0$ is reduced. For the case 
of super-spreaders, the curves always start with a decreasing slope 
and later change the peak but with lower total number of infected individuals. 
This makes these classes particularly special and might have a huge effect 
in the progression of the other infected classes. In the particular 
case of Portugal, it is remarkable that the peak of infected super-spreaders 
is less than one.  


\clearpage
\end{paracol}
\nointerlineskip
\begin{paracol}{2}

\switchcolumn


\clearpage
\end{paracol}
\nointerlineskip

\begin{paracol}{2}

\switchcolumn


\clearpage
\end{paracol}
\nointerlineskip
\begin{paracol}{2}

\switchcolumn
\vspace{-6pt}


\section{Conclusions}
\label{sec:conc}

In this paper, we have analyzed in detail a fractional-order COVID-19 model, previously 
used for fitting the number of confirmed cases of infections from the region 
of Galicia, Spain, and Portugal \cite{fays2}. An analytical study on the stability 
of the steady state was conducted and numerical simulations investigated 
for all infected compartments of the population. The model qualitative 
analysis reconfirm stability of the steady state for both classical 
and fractional-order models, whenever the threshold condition $R_0<1$ holds.
Global stability is then deduced by a fractional-order version of 
LaSalle's invariant set theorem. Moreover, our numerical solutions show that 
the time speed evolution for the dynamical model to reach the steady state 
is affected by the order $\alpha$ of index memory. It remains open the question 
of how to prove global stability in the case when $R_0>1$. For all numerical 
simulations we have done, the equilibrium was always stable.


\vspace{6pt} 

\authorcontributions{Conceptualization, F.N.; methodology, F.N. and D.F.M.T.; 
software, F.N.; validation, F.N. and D.F.M.T.; formal analysis, F.N. and D.F.M.T.; 
investigation, F.N. and D.F.M.T.; writing---original draft preparation, F.N. and D.F.M.T.;
writing---review and editing, F.N. and D.F.M.T.; supervision, D.F.M.T. 
All authors have read and agreed to the published version of the~manuscript.}

\funding{This research was funded by The Portuguese Foundation for Science 
and Technology (FCT) grant number UIDB/04106/2020 (CIDMA).
Fa\"{\i}\c{c}al Nda\"{\i}rou was also supported by FCT through
the PhD fellowship PD/BD/150273/2019.}

\institutionalreview{Not applicable.}

\informedconsent{Not applicable.}

\dataavailability{Not applicable.} 

\acknowledgments{The authors are grateful to two reviewers for their comments and suggestions.}

\conflictsofinterest{The authors declare no conflicts of interest. 
The funders had no role in the design of the study; in the collection, 
analyses, or interpretation of data; in the writing of the manuscript, 
or in the decision to publish the~results.} 


\end{paracol}


\reftitle{References}

 

\begin{thebibliography}{999}

\bibitem[Oukouomi~Noutchie \em{et~al.}(2014)Oukouomi~Noutchie, Atangana, and Goufo]{atangana2014}
Oukouomi~Noutchie, S.C.; Atangana, A.; Goufo, E.F.D.
\newblock On the Mathematical Analysis of Ebola Hemorrhagic Fever: Deathly
Infection Disease in West {A}frican Countries.
\newblock {\em BioMed Res. Int.} {\bf 2014}, {\em 2014},~261383.

\bibitem[Yanover \em{et~al.}(2016)Yanover, Ngwa, and Teboh-Ewungkem]{Ngwa}
Yanover, C.; Ngwa, G.A.; Teboh-Ewungkem, M.I.
\newblock A Mathematical Model with Quarantine States for the Dynamics of Ebola
Virus Disease in Human Populations.
\newblock {\em Comput. Math. Methods Med.} {\bf 2016}, {\em 2016},~9352725.

\bibitem[Yuan \em{et~al.}(2015)Yuan, Rachah, and Torres]{MR3349757}
Yuan, S.; Rachah, A.; Torres, D.F.M.
\newblock Mathematical Modelling, Simulation, and Optimal Control of the 2014
Ebola Outbreak in West Africa.
\newblock {\em Discret. Dyn. Nat. Soc.} {\bf 2015}, {\em 2015},~842792.
\newblock {\tt arXiv:1503.07396}

\bibitem[Djordjevic \em{et~al.}(2018)Djordjevic, Silva, and Torres]{MR3808514}
Djordjevic, J.; Silva, C.J.; Torres, D.F.M.
\newblock A stochastic SICA epidemic model for HIV transmission.
\newblock {\em Appl. Math. Lett.} {\bf 2018}, {\em 84},~168--175.
\newblock doi:10.1016/j.aml.2018.05.005.
\newblock {\tt arXiv:1805.01425}

\bibitem[Duwal \em{et~al.}(2015)Duwal, Winkelmann, Sch\"{u}tte, and von Kleist]{OC:HIV:PLoSCompBio:2015}
Duwal, S.; Winkelmann, S.; Sch\"{u}tte, C.; von Kleist, M.
\newblock Optimal treatment strategies in the context of 'Treatment for
Prevention' against HIV-1 in resource-poor settings.
\newblock {\em PLoS Comput. Biol.} {\bf 2015}, {\em 11},~e1004200.

\bibitem[Silva and Torres(2015)]{SilvaTorres:TBHIV:2015}
Silva, C.J.; Torres, D.F.M.
\newblock A TB-HIV/AIDS coinfection model and optimal control treatment.
\newblock {\em Discret. Contin. Dyn. Syst.  A} {\bf 2015}, {\em 35},~4639--4663.
\newblock {\tt arXiv:1501.03322}

\bibitem[Hu \em{et~al.}(2012)Hu, Ma, and Ruan]{HU201212}
Hu, Z.; Ma, W.; Ruan, S.
\newblock Analysis of SIR epidemic models with nonlinear incidence rate and treatment.
\newblock {\em Math. Biosci.} {\bf 2012}, {\em 238},~12--20.
\newblock doi:10.1016/j.mbs.2012.03.010.

\bibitem[Hethcote and van~den Driessche(1991)]{PMID:2061695}
Hethcote, H.W.; van~den Driessche, P.
\newblock Some epidemiological models with nonlinear incidence.
\newblock {\em J. Math. Biol.} {\bf 1991}, {\em 29},~271--287.
\newblock doi:{\changeurlcolor{black}\href{https://doi.org/10.1007/bf00160539}{\detokenize{10.1007/bf00160539}}}.

\bibitem[Moreno \em{et~al.}(2017)Moreno, Espinoza, Bichara, Holechek, and Castillo-Chavez]{Moreno:2017}
Moreno, V.M.; Espinoza, B.; Bichara, D.; Holechek, S.A.; Castillo-Chavez, C.
\newblock Role of short-term dispersal on the dynamics of Zika virus in an extreme idealized environment.
\newblock {\em Infect. Dis. Model.} {\bf 2017}, {\em 2},~21--34.

\bibitem[Nda{\"\i}rou \em{et~al.}(2018)Nda{\"\i}rou, Area, Nieto, Silva, and Torres]{fays4}
Nda{\"\i}rou, F.; Area, I.; Nieto, J.J.; Silva, C.J.; Torres, D.F.M.
\newblock Mathematical modeling of Zika disease in pregnant women and newborns
with microcephaly in Brazil.
\newblock {\em Math. Methods Appl. Sci.} {\bf 2018}, {\em 41},~8929--8941.
\newblock {\tt arXiv:1711.05630}

\bibitem[Nda{\"\i}rou \em{et~al.}(2020)Nda{\"\i}rou, Area, and Torres]{fays3}
Nda{\"\i}rou, F.; Area, I.; Torres, D.F.M.
\newblock Mathematical Modeling of Japanese Encephalitis under Aquatic Environmental Effects.
\newblock {\em Mathematics} {\bf 2020}, {\em 8}, 1880.
\newblock doi:{\changeurlcolor{black}\href{https://doi.org/10.3390/math8111880}{\detokenize{10.3390/math8111880}}}.
\newblock {\tt arXiv:2010.09418}

\bibitem[Rodrigues \em{et~al.}(2016)Rodrigues, Monteiro, and Torres]{MR3557143}
Rodrigues, H.S.; Monteiro, M.T.T.; Torres, D.F.M.
\newblock Seasonality effects on dengue: Basic reproduction number, sensitivity analysis and optimal control.
\newblock {\em Math. Methods Appl. Sci.} {\bf 2016}, {\em 39},~4671--4679.
\newblock {\tt arXiv:1409.3928}

\bibitem[Kermack and McKendrick(1927)]{Kendrick1}
Kermack, W.O.; McKendrick, A.G.
\newblock Contributions to the mathematical theory of epidemics I.
\newblock {\em Proc. R. Soc. Lond.} {\bf 1927}, {\em 115},~700--721.

\bibitem[Kermack and McKendrick(1932)]{Kendrick2}
Kermack, W.O.; McKendrick, A.G.
\newblock Contributions to the mathematical theory of epidemics II.
\newblock {\em Proc. R. Soc. Lond.} {\bf 1932}, {\em 138},~55--83.

\bibitem[Kermack and McKendrick(1933)]{Kendrick3}
Kermack, W.O.; McKendrick, A.G.
\newblock Contributions to the mathematical theory of epidemics III.
\newblock {\em Proc. R. Soc. Lond.} {\bf 1933}, {\em 141},~94--112.

\bibitem[Area \em{et~al.}(2017)Area, Losada, Nda{\"\i}rou, Nieto, and Tcheutia]{Area:in:press}
Area, I.; Losada, J.; Nda{\"\i}rou, F.; Nieto, J.J.; Tcheutia, D.D.
\newblock Mathematical modeling of 2014 Ebola outbreak.
\newblock {\em Math. Methods Appl. Sci.} {\bf 2017}, {\em 40},~6114--6122.

\bibitem[Jordan and Smith(1999)]{smith}
Jordan, Y.A.; Smith, P.
\newblock {\em Nonlinear Ordinary Differential Equations. An Introduction to
Dynamical Systems}; Oxford University Press Inc.: New York, NY, USA, 1999.

\bibitem[Melnik and Korobeinikov(2013)]{melnikkk}
Melnik, A.V.; Korobeinikov, A.
\newblock Lyapunov functions and global stability for SIR and SEIR models with age-dependent susceptibility.
\newblock {\em Math. Biosci. Eng.} {\bf 2013}, {\em 10},~369--378.

\bibitem[Strogatz(1994)]{steven}
Strogatz, S.H.
\newblock {\em Nonlinear Dynamics and Chaos: With Applications to Physics,
Biology, Chemistry, and Engineering}; Addison-Wesley Pub: Reading, MA, USA, 1994.

\bibitem[Teschl(2012)]{teschl}
Teschl, G.
\newblock {\em Ordinary Differential Equations and Dynamical Systems}; 
American Mathematical Society: Providence, RI, USA, 2012.

\bibitem[Al-{S}ulami \em{et~al.}(2014)Al-{S}ulami, El-{S}hahed, and Nieto]{sulami}
Al-{S}ulami, H.; El-{S}hahed, M.; Nieto, J.J.
\newblock On fractional order dengue epidemic model.
\newblock {\em Math. Probl. Eng.} {\bf 2014}, doi:10.1155/2014/456537.

\bibitem[Rosa and Torres(2019)]{MR3999702}
Rosa, S.; Torres, D.F.M.
\newblock Optimal control and sensitivity analysis of a fractional order {TB} model.
\newblock {\em Stat. Optim. Inf. Comput.} {\bf 2019}, {\em 7},~617--625.
\newblock doi:{\changeurlcolor{black}\href{https://doi.org/10.19139/soic.v7i3.836}{\detokenize{10.19139/soic.v7i3.836}}}.
\newblock {\tt arXiv:1812.04507}

\bibitem[Nda{\"\i}rou \em{et~al.}(2021)Nda{\"\i}rou, Area, Nieto, Silva, and Torres]{fays2}
Nda{\"\i}rou, F.; Area, I.; Nieto, J.J.; Silva, C.J.; Torres, D.F.M.
\newblock Fractional model of COVID-19 applied to Galicia, Spain and Portugal.
\newblock {\em Chaos Solitons Fractals} {\bf 2021}, {\em 144},~110652.
\newblock doi:10.1016/j.chaos.2021.110652.
\newblock {\tt arXiv:2101.01287}

\bibitem[Du \em{et~al.}(2013)Du, Wang, and Hu]{memory}
Du, M.; Wang, Z.; Hu, H.
\newblock Measuring memory with the order of fractional derivative.
\newblock {\em Sci. Rep.} {\bf 2013}, {\em 3}, 3431.

\bibitem[Noeiaghdam \em{et~al.}(2021)Noeiaghdam, Micula, and Nieto]{math9121321}
Noeiaghdam, S.; Micula, S.; Nieto, J.J.
\newblock A Novel Technique to Control the Accuracy of a Nonlinear Fractional
Order Model of {COVID-19}: Application of the {CESTAC} Method and the {CADNA} Library.
\newblock {\em Mathematics} {\bf 2021}, {\em 9}, 1321.
\newblock doi:{\changeurlcolor{black}\href{https://doi.org/10.3390/math9121321}{\detokenize{10.3390/math9121321}}}.

\bibitem[Rosa and Torres(2018)]{MR3872489}
Rosa, S.; Torres, D.F.M.
\newblock Optimal control of a fractional order epidemic model with application
to human respiratory syncytial virus infection.
\newblock {\em Chaos Solitons Fractals} {\bf 2018}, {\em 117},~142--149.
\newblock doi:{\changeurlcolor{black}\href{https://doi.org/10.1016/j.chaos.2018.10.021}{\detokenize{10.1016/j.chaos.2018.10.021}}}.
\newblock {\tt arXiv:1810.06900}

\bibitem[Boukhouima \em{et~al.}(2021)Boukhouima, Lotfi, Mahrouf, Rosa, Torres, and Yousfi]{Boukhouima2021}
Boukhouima, A.; Lotfi, E.M.; Mahrouf, M.; Rosa, S.; Torres, D.F.M.; Yousfi, N.
\newblock Stability analysis and optimal control of a fractional {HIV}-{AIDS} 
epidemic model with memory and general incidence rate.
\newblock {\em Eur. Phys. J. Plus} {\bf 2021}, {\em 136}. 
\newblock
doi:{\changeurlcolor{black}\href{https://doi.org/10.1140/epjp/s13360-020-01013-3}{\detokenize{10.1140/epjp/s13360-020-01013-3}}}.
\newblock {\tt arXiv:2012.04819}

\bibitem[Samko \em{et~al.}(1993)Samko, Kilbas, and Marichev]{samko1993fractional}
Samko, S.G.; Kilbas, A.A.; Marichev, O.I.
\newblock {\em Fractional Integrals and Derivatives}; 
Gordon and Breach Science Publishers, Yverdon, 1993. 

\bibitem[Kilbas \em{et~al.}(2006)Kilbas, Srivastava, and Trujillo]{MR2218073}
Kilbas, A.A.; Srivastava, H.M.; Trujillo, J.J.
\newblock {\em Theory and Applications of Fractional Differential Equations};
{North-Holland Mathematics Studies}; Elsevier Science B.V.: 
Amsterdam, The Netherlands, 2006.

\bibitem[Podlubny(1999)]{MR1658022}
Podlubny, I.
\newblock {\em Fractional Differential Equations}; 
Academic Press, Inc.: San Diego, CA,  USA, 1999.

\bibitem[Sidi~Ammi \em{et~al.}(2021)Sidi~Ammi, Tahiri, and Torres]{MR4232864}
Sidi~Ammi, M.R.; Tahiri, M.; Torres, D.F.M.
\newblock Global {S}tability of a {C}aputo {F}ractional {SIRS} {M}odel 
with {G}eneral {I}ncidence {R}ate.
\newblock {\em Math. Comput. Sci.} {\bf 2021}, {\em 15},~91--105.
\newblock
doi:{\changeurlcolor{black}\href{https://doi.org/10.1007/s11786-020-00467-z}{\detokenize{10.1007/s11786-020-00467-z}}}.
\newblock {\tt arXiv:2002.02560}

\bibitem[Lin(2007)]{lin}
Lin, W.
\newblock Global existence theory and chaos control of fractional differential equations.
\newblock {\em J. Math. Anal. Appl.} {\bf 2007}, {\em 332},~709--726.

\bibitem[Almeida \em{et~al.}(2019)Almeida, Brito~da Cruz, Martins, and Monteiro]{ricardo}
Almeida, R.; Brito~da Cruz, A.M.C.; Martins, N.; Monteiro, M.T.T.
\newblock An epidemiological MSEIR model described by the Caputo fractional derivative.
\newblock {\em Int. J. Dyn. Control} {\bf 2019}, {\em 7},~776--784.

\bibitem[Nda{\"\i}rou \em{et~al.}(2020)Nda{\"\i}rou, Area, Nieto, and Torres]{fays1}
Nda{\"\i}rou, F.; Area, I.; Nieto, J.J.; Torres, D.F.M.
\newblock {Mathematical modeling of COVID-19 transmission dynamics with a case study of Wuhan}.
\newblock {\em Chaos Solitons Fractals} {\bf 2020}, {\em 135},~109846.
\newblock {\tt arXiv:2004.10885}

\bibitem[Nda{\"\i}rou \em{et~al.}(2020)Nda{\"\i}rou, Area, Bader, Nieto, and Torres]{MyID:460}
Nda{\"\i}rou, F.; Area, I.; Bader, G.; Nieto, J.J.; Torres, D.F.M.
\newblock {Corrigendum to 'Mathematical Modeling of COVID-19 Transmission 
Dynamics with a Case Study of Wuhan' [Chaos Solitons Fractals 135 (2020), 109846]}.
\newblock {\em Chaos Solitons Fractals} {\bf 2020}, {\em 141},~110311.

\bibitem[van~den Driessche and Watmough(2002)]{van:den:Driessche:2002}
Van~den Driessche, P.; Watmough, J.
\newblock Reproduction numbers and sub-threshold endemic equilibria 
for compartmental models of disease transmission.
\newblock {\em Math. Biosci.} {\bf 2002}, {\em 180},~29--48. 
doi:10.1016/S0025-5564(02)00108-6.

\bibitem[Diekmann \em{et~al.}(1990)Diekmann, Heesterbeek, and Metz]{diekmann}
Diekmann, O.; Heesterbeek, J.A.P.; Metz, J.A.J.
\newblock On the definition and the computation of the basic reproduction ratio
${R}_{0}$ in models for infectious diseases in heterogeneous populations.
\newblock {\em J. Math. Biol.} {\bf 1990}, {\em 28},~365--382.

\bibitem[Huo \em{et~al.}(2015)Huo, Zhao, and Zhu]{MR3384337}
Huo, J.; Zhao, H.; Zhu, L.
\newblock The effect of vaccines on backward bifurcation in a fractional order {HIV} model.
\newblock {\em Nonlinear Anal. Real World Appl.} {\bf 2015}, {\em 26},~289--305.
\newblock
doi:{\changeurlcolor{black}\href{https://doi.org/10.1016/j.nonrwa.2015.05.014}{\detokenize{10.1016/j.nonrwa.2015.05.014}}}.

\end{thebibliography}
\end{document}